\documentclass[12pt,a4paper,twoside,final,notitlepage, leqno]{article}
\usepackage[english]{babel}
\usepackage[T1]{fontenc}  
\usepackage{epsfig, graphicx, amssymb}
\usepackage{amsmath,amsthm,epsfig,amsfonts}  
\usepackage{float}
\usepackage{color}
\def \smb {{\scriptstyle \bullet }}
 
\newcommand{\moneq}{\vspace*{-7pt} \begin{equation} \displaystyle } 
\newcommand{\moneqstar}{\vspace*{-6pt} \begin{equation*} \displaystyle } 
\newcommand{\monendstar}{\vspace*{-6pt} \end{equation*}   }
\newcommand{\monend}{\vspace*{-7pt} \end{equation}   }

\def\eq{\mathop{\rm {eq}}\nolimits}

\def\section*#1{}

\usepackage{fancyhdr}
\renewcommand{\headrulewidth}{0pt}

\begin{document} 

\fancypagestyle{plain}{ \fancyfoot{} \renewcommand{\footrulewidth}{0pt}}
\fancypagestyle{plain}{ \fancyhead{} \renewcommand{\headrulewidth}{0pt}}

~


\centerline {\bf \Large   Recovering the full Navier Stokes equations   }

 \bigskip 

\centerline {\bf \Large  with lattice Boltzmann schemes  }

 \bigskip  \bigskip \bigskip

\centerline { \large    Fran\c{c}ois Dubois$^{ab}$, Benjamin Graille$^{a}$, Pierre Lallemand$^{c}$ }

\smallskip  \bigskip 

\centerline { \it  \small   
$^a$   Dpt. of Mathematics, University Paris-Sud,  B\^at. 425, F-91405  Orsay, France.} 

\centerline { \it  \small   
$^b$    Conservatoire National des Arts et M\'etiers, LMSSC laboratory,  F-75003 Paris, France.} 

\centerline { \it  \small  $^c$   Beijing Computational Science Research Center, 
Haidian District, Beijing 100094,  China.}


\bigskip  

\centerline { 13 July 2016  
{\footnote {\rm  \small $\,$ Published in November 2016 by the 
in American Institute of Physics Proceedings, volume 1786, 040003 (2016), 
doi.org/10.1063/1.4967541.  Edition 20 November 2017.}}}

 \bigskip 


\bigskip  
\noindent {\bf  Abstract} 

\noindent 
We  consider  multi relaxation times  
lattice Boltzmann scheme 
with two particle distributions 
for the thermal  Navier Stokes equations formulated with 
conservation of mass and momentum and  dissipation of volumic entropy.
Linear stability is  taken into consideration 
to determine a coupling between two coefficients of dissipation.
We present interesting  numerical results 
for one-dimensional strong nonlinear acoustic waves with shocks.

\bigskip \bigskip   \noindent {\bf      1) \quad  Introduction  }    

\noindent 
In this contribution, our program is 
to use the lattice Boltzmann schemes with multiresolution relaxation times \cite{DdH92}  
for the approximation of the full compressible Navier Stokes equations. 
Interesting  results have been obtained previously in \cite{LD15,LL03} 
when using a single particle distribution. 
A popular idea proposed in  \cite{ACS93,KP94,MNA93,SC93}  is the
use of several particle distributions. One particle distribution is devoted to the conservation of mass
and momentum and an other one to the conservation of energy.
%
This framework has been used in the context 
of Bhatnagar Gross Krook \cite{BGK54} approximation \cite{GZS07,KSC13,LHWT07,LLHT12}.
The adaptation to the  multiresolution relaxation times approach is not straightforward  \cite{Du16} 
and  source terms must be added to a pure
``collide-stream'' algorithm in order to capture this nonlinear dissipation term.

\smallskip \noindent
In this contribution, we begin with elementary one-dimensional fluid flow 
with the D1Q3 multi relaxation times lattice Boltzmann scheme. 
The Taylor expansion method allows us to determine the underlying partial differential equations. 
We adapt this study to the case of only one advection-diffusion equation. 
The coupling of these two schemes for the approximation
of the Navier Stokes equations conducts to difficulties.
 This defect has its origin in precise algebraic properties of  lattice
Boltzmann schemes and associated moments.  We then admit to consider 
lattice Boltzmann schemes with source terms. We focus on a formulation of the 
 Navier Stokes equations that takes into consideration the production of entropy. 
The coupling of two D1Q3 schemes is proposed and studied and first numerical results
are presented.

\bigskip \bigskip   \noindent {\bf    2) \quad  
Fluid flow with a multi relaxation times D1Q3 lattice Boltzmann scheme}  

\noindent 
We describe in this section a very simple but fundamental lattice Boltzmann scheme. 
We restric ourselves to one space dimension. The  mesh 
is parameterized by a space step  $\Delta x$. 
The nodes are located at the vertices $x= j \, \Delta x$, where $j$ is an integer. 
We suppose given three discrete velocities:   $ v_0 = 0 $, $ v_+ = +1 $ and  $ v_- = -1 $. 
A first description of the degrees of freedom is associated to the density 
$ f \equiv  ( f_0 ,\, f_+ ,\, f_- ) $  of particles.
In the  D1Q3 case, we have three  kind of particles. 
The  motionless particles $f_0(x,t)$, 
the  particles $f_+(x,t)$ going from  $x$ to $ x + \Delta x $ during one time step $ \Delta t $ 
and the particles $f_-(x,t)$ going from  $x$ to $ x - \Delta x $ during the time step. 
We adopt the so-called ``acoustic scale'' and 
the numerical velocity $ \lambda \equiv {{\Delta x}\over{\Delta t}}  $ 
is supposed fixed for the entire contribution. 
From this distribution of particle densities, we construct two first moments: 
mass density  
\moneqstar 
 \rho   \equiv   \rho_0  \sum_{j=0 ,\, + ,\, -}  f_j 
 =  \rho_0  \, (  f_+ + f_- + \, f_0 )   
\monendstar 
and   momentum 
\moneqstar 
 J    \equiv    \rho_0 \, \lambda \, \sum_{j=0 ,\, + ,\, -}  v_j \,\, f_j 
 =  \rho_0 \, \lambda \, (f_+ - f_- )  \, . 
\monendstar 
The  third momentum is the so-called ``energy'' defined by 
\moneqstar 
 e   \equiv   \rho_0 \, \lambda^2  \, (  f_+ + f_- - 2 \, f_0 ) \, . 
\monendstar 
The second description of the degrees of freedom is due to the moments 
\moneqstar 
 m \equiv ( \rho ,\, J ,\, e )  \, . 
\monendstar 
An invertible matrix $M$ is therefore defined between the particles and the moments: 
\moneq
 m \,\equiv\, \left( \begin{array}{c}  \rho \\  J \\ e \end{array} \right) 
\,=\, M_{\rm D1Q3}  \, \,  
\left( \begin{array}{c}  f_0 \\ f_+ \\ f_- \\  \end{array} \right)  
\,\equiv \,  M_{\rm D1Q3}  \, f \,, \quad 
  M_{\rm D1Q3} \,=\, \rho_0 \, \left( 
\begin{array}{ccc}
   1 &   1 & 1  \\ 
   0 &   \lambda & -\lambda  \\ 
   -2 \, \lambda^2  & \lambda^2 &  \lambda^2  
\end{array} \right)  \, . \label{M-d1q3} 
\monend 
The evolution algorithm is composed by two steps: relaxation and advection.
In the relaxation step, 
density and momentum remain at  equilibrium and do not change: 
\moneqstar 
 \rho^* = \rho \,,\,\,  J^* = J \, . 
\monendstar 
At the contrary, 
the equilibrium energy is a given function of the two moments at equilibrium: 
\moneqstar 
 e^{\rm eq} =  e^{\rm eq} (\rho , \, J)  \, . 
\monendstar 
%
The energy after relaxation is obtained by a simple 
evolution: 
\moneqstar 
  e^* =  e + s_e \, (e^{\rm eq} - e)  \, . 
\monendstar 
 The relaxation 
parameter $s_e$ must be chosen satisfying  $ 0 < s_e < 2 $ \cite{LL00}. 
The particle distribution $f^*$ after relaxation is defined by 
\moneqstar 
f^* =  M_{\rm D1Q3} ^{-1} \, (\rho ,\,  J ,\,  e^* )^{\rm t}  \, . 
\monendstar 
The advection step is a free displacement of the particles 
during one time step: 
\moneqstar 
 f_0 (x, \, t+\Delta t) =  f_0^* (x, \, t) \,,\,\, 
   f_+ (x, \, t+\Delta t) =  f_+^* (x - \Delta x ,\, t) \,,\,\, 
  f_- (x, \, t+\Delta t) =  f_-^* (x + \Delta x ,\, t) \, . 
\monendstar

\fancyhead[EC]{\sc{ Fran\c c ois Dubois,  Benjamin Graille and  Pierre Lallemand }} 
\fancyhead[OC]{\sc{ Recovering the full Navier Stokes equations with lattice Boltzmann }} 
\fancyfoot[C]{\oldstylenums{\thepage}}

\bigskip \bigskip   \noindent {\bf     3) \quad  
 Taylor expansion method  and associated algebraic tools }  

\noindent 
We conduct the analysis of the  multi relaxation times
with the Taylor expansion method \cite{Du08}. 
This approach is a kind of numerical Chapman-Enskog expansion 
with the space step $\Delta x $ as a small parameter. 
Equivalent partial differential equations emerge from this analysis for the 
conserved variables. We have for the previous example up to order 2:
\moneq
\partial_t \rho + \partial_x J = {\rm O}(\Delta x^2) \,,\quad  
\partial_t J + 
\partial_x  \Big( {2\over3} \lambda^2 \, \rho  + {1\over3} e^{\rm eq}  \Big) 
- {1\over3}  \,  \sigma_e \, \Delta t \, \partial_x \theta_e = {\rm O}(\Delta x^2)  \, , 
 \label{edp-d1q3}  \monend 
with a coefficient 
$ \displaystyle \sigma_e = {{1}\over{s_e}} - {1\over2} $ introduced by 
H\'enon \cite{He87}  
and 
\moneqstar 
  \theta_e  \equiv \partial_t  e^{\rm eq}  +  \lambda^2 \, \partial_x J  
\simeq  3 \, \lambda^2  \rho \, \partial_x  u 
\monendstar  
for the previous fluid D1Q3 scheme. 
We can compare these equations (\ref{edp-d1q3}) with the 
unidimensional  ``isentropic'' Navier Stokes equations 
\moneqstar 
 \partial_t \rho + \partial_x J =  0  
\monendstar  
 for mass conservation and 
\moneqstar 
 \partial_t J + \partial_x \big( \rho \, u^2 + p \big) 
- \partial_x \big( \rho \, \nu \, \partial_x u \big) =  0  
\monendstar  
for momentum conservation. 
Then the 
energy at equilibrium follows the relation
\moneq  
e^{\eq} = 3 \, \big( \rho \, u^2 + p \big) - 2 \, \lambda^2 \, \rho \, . 
 \label{energie-equilibre}  \monend 
The relaxation parameter is defined in (\ref{henon}) and 
is related to  kinematic viscosity: 
\moneqstar 
   \nu_0  = \lambda  \, \Delta x \, \sigma_e   \, . 
\monendstar  
In all generality, if the moments $m$ are related to the particle distribution $f$ 
{\it via} a relation of the type 
$  m =  M \, f   $ 
with a invertible matrix $M$, we introduce  \cite{Du08} the 
momentum-velocity tensor  $ \Lambda $ according to 
\moneq   
\Lambda_{k \ell} \equiv  \lambda \, \sum_{j}  M_{k  j}  \, v_j \, (M^{-1})_{j  \ell} \, . 
 \label{Lambda-general}  \monend 
For the previous D1Q3 model,  with the ordering $ \rho ,\,  J ,\,  e  $  of the moments $m$, 
we have  
\moneq
\Lambda_{\rm D1Q3}  =    \left( \begin{array}{ccc} 
   0  &   1 & 0  \\ 
  {2\over3} \lambda^2 &  0  & {1\over3}  \\
  0 & \lambda^2 & 0    \end{array} \right) \, . 
 \label{tenseur-Lambda}  \monend 
The vector of conserved moments is denoted by  $ W $ 
and the defect of conservation  $ \theta_\ell $ is defined according  to 
\moneq
\theta_\ell \equiv \partial_t  m_\ell^{\rm eq} +  \sum_p  \Lambda_{\ell p} \,  \partial_x m_p^{\rm eq}
 \label{thetas}  \monend 
for the non-conserved moments $ m_\ell $. Then the $ k^{\rm th} $ conserved moment $ W_k $ satisfies 
asymptotically the  following  second order  partial differential equation 
\moneq 
\partial_t   W_k  +  \sum_\ell \Lambda_{k \ell}  \,\,  \partial_x m_\ell^{\rm eq}
- \Delta t \,\, \sum_\ell \, \sigma_\ell \,\, \Lambda_{k \ell} 
\,\, \partial_x  \theta_\ell   =  {\rm O}(\Delta x^2) \, . 
 \label{edp-ordre-2}  \monend 
with the H\'enon coefficient \cite{He87} $ \sigma_\ell $ defined in all generality 
according to 
\moneq 
\sigma_\ell = {{1}\over{s_\ell}} - {1\over2} \, . 
 \label{henon}  \monend 
The generalization to two and three dimensions is straightforward.
 We observe also that 
the equations (\ref{edp-ordre-2}) are always under a conservative form.

\bigskip \bigskip   \noindent {\bf   4) \quad  
 Advection-diffusion with a multi relaxation times D1Q3 

\quad $\,\,\,\,\,\,$ lattice Boltzmann scheme }  

\noindent 
We modify the notations of the previous D1Q3 ``MRT'' scheme, replacing the 
particle distribution $f$ by the notation $g$:
$g_0(x,t)$  is the  motionless particles at   $x$ during the time step  $ \Delta t $, 
$g_+(x,t)$  the  density of particles going from  $x$ to $ x + \Delta x $ during  one  time step 
and  $g_-(x,t)$  the density of particles going from  $x$ to $ x - \Delta x $ during $ \Delta t $.
The numerical velocity $  \lambda  $ remains fixed.  
We introduce three  moments: a  conserved variable  $ \zeta$, momentum $\psi$ and associated energy  
$\varepsilon$ in a way analogous to (\ref{M-d1q3}):
\moneqstar 
\zeta \equiv   \zeta_0  \, (   g_0 + g_+ + g_-  ) \,,\,\,
  \psi    \equiv     \zeta_0 \, \lambda \, (g_+ -  g_- )  \,,\,\, 
 \varepsilon  \equiv   \zeta_0  \, \lambda^2  \, (  g_+ + g_- - 2 \, g_0 )  \,. 
\monendstar 
An invertible matrix  $\widetilde{M}$ is defined by the condition 
\moneqstar 
 (  \zeta ,\, \psi ,\, \varepsilon )^{\rm t}  = \widetilde{M}_{\rm D1Q3}  \, 
(  g_0 ,\,  g_+  ,\,g_-  )^{\rm t} \,. 
\monendstar 
It is analogous to the matrix introduced in (\ref{M-d1q3}) except that the scaling factor 
$ \rho_0 $  is replaced by $ \zeta_0 $. 
For this model, only one moment is conserved and  remains at   equilibrium: 
$ \zeta^* = \zeta $.
The momentum  $ \psi  $ and associated  energy $   \varepsilon  $ 
at equilibrium  are given functions  of the scalar~$ \zeta $:  
\moneqstar 
  \psi^{\rm eq} =  \psi^{\rm eq} (\zeta) \,,\,\, 
  \varepsilon^{\rm eq} =   \varepsilon^{\rm eq} (\zeta)     \,. 
\monendstar 
The momentum and energy after relaxation $  \psi^* $ and $   \varepsilon^* $
are  obtained  by simple evolutions: 
\moneqstar 
  \psi^* =  \psi  + s_\psi \, (\psi^{\rm eq} - \psi)   \,\,\,
   \varepsilon^* =  \varepsilon  + s_\varepsilon \, (\varepsilon^{\rm eq} - \varepsilon)  
\monendstar 
 with  $ 0 < s_\psi < 2 $ and $   0 < s_\varepsilon < 2 $.  
The particle distribution $ g^* $ after relaxation is  obtained from the moments:
\moneqstar 
 ( g_0^* ,\,  g_+^* ,\,  g_-^* )^{\rm t} = ( \widetilde{M}_{\rm D1Q3} ) ^{-1} \, 
 (  \zeta  ,\,  \psi^*  ,\,  \varepsilon^* )^{\rm t}  \,. 
\monendstar 
The advection step is identical to the corresponding one for the fluid scheme:
\moneqstar  \left\{ \begin{array} {l} 
\displaystyle 
  g_0 (x, \, t+\Delta t) =  g_0^* (x, \, t)   \\   
  g_+ (x, \, t+\Delta t) =  g_+^* (x - \Delta x ,\, t)  \\
  g_- (x, \, t+\Delta t) =  g_-^* (x + \Delta x ,\, t) \, . 
\end{array}  \right. \monendstar 
%
We can analyze formally this scheme with the  Taylor expansion method. 
With the matrix $ \Lambda $ introduced in (\ref{tenseur-Lambda}), 
the equivalent partial differential equation up to order 2 takes the form
\moneq 
 \partial_t \zeta + \partial_x \psi^{\rm eq}  
- \sigma_\psi \, \Delta t \, \partial_x \theta_\psi = {\rm O}(\Delta x^2) 
 \label{edp-thermic-ordre-2}  \monend 
with  
\moneqstar 
 \sigma_\psi = {{1}\over{s_\psi}} - {1\over2} 
\monendstar 
and  
\moneqstar 
  \theta_\psi  \equiv \partial_t  \psi^{\rm eq}  + 
\partial_x  \Big( {2\over3} \lambda^2 \, \rho  + {1\over3} \varepsilon^{\rm eq}   \Big) \, . 
\monendstar 
This ``thermal'' lattice Boltzmann model is well adapted for the simulation 
of an advection-diffusion equation of the type 
\moneqstar 
 \partial_t \zeta + \partial_x  ( u_0 \, \zeta  )  
- \partial_x  ( \kappa \, \partial_x \zeta  ) =  0   \, . 
\monendstar 
The second order differential equation  (\ref{edp-thermic-ordre-2}) 
simulated by the thermal lattice Boltzmann scheme  is a good approximation to  first  
order of (\ref{edp-thermic-ordre-2})
if the momentum at equilibrium  $ \psi^{\eq} $ is given by 
\moneqstar  
\psi^{\eq} = u_0 \,\, \zeta  \, . 
\monendstar 
The defect of equilibrium  can be expanded: 
\moneqstar  
 \theta_\psi \equiv   \Big(  {2\over3} \lambda^2 - u_0^2  \Big) \, \partial_x \zeta 
+  {1\over3} \, \partial_x \varepsilon^{\rm eq} + {\rm O}(\Delta x) 
\monendstar 
and the energy  at equilibrium is  proportional to the conserved variable: 
\moneqstar  
  \varepsilon^{\eq} = \alpha \, \lambda^2  \, \zeta  \, . 
\monendstar 
The identification of second order terms induces:
\moneqstar  
    \Big(  {{2+\alpha}\over3} \lambda^2 - u_0^2   \Big) \,\,  
 \Big(  {{1}\over{s_\psi}} - {1\over2} \Big) \,\, \Delta t  = \kappa 
\monendstar 
with a typical constraint for numerical stability: 
  $     -2  < \alpha < 1   $.  
We take simply $ \, s_\varepsilon = 1.5 \,$ for the third moment.

\bigskip \bigskip   \noindent {\bf    5) \quad  
 Unidimensional Navier Stokes equations with energy conservation }   

\noindent 
We try now to simulate the one dimensional Navier Stokes equations with the coupling the two previous 
models. The conserved variables
are the volumic mass  $ \rho  $, the momentum  $ J \, \equiv \, \rho  \, u \,, $
and the  volumic total  energy 
\moneqstar  
 \rho  \, E \equiv  \rho  \,  \Big( i + {1\over2} \, u^2  \Big) 
\monendstar 
obtained by adding the specific kinetic energy $ {1\over2} \, u^2 $ to the internal one $i$.
We suppose that the  equation of state   
is a polytropic perfect gas with a ratio
\moneqstar  
  \gamma \equiv {{c_p}\over{c_v}} 
\monendstar 
 of specific heats: 
\moneqstar  
 p =  (\gamma - 1)  \, \rho \, i  =  \rho \, r \, T  \, . 
\monendstar 
The expression of the sound velocity  $c$ is classical: 
\moneqstar  
 c^2  =  \gamma \, {{p}\over{\rho}} =   \gamma \, (\gamma - 1)  \, i  
 =   (\gamma - 1)  \, c_p \, \, T  \, . 
\monendstar 
The kinematic viscosity and thermal conductivity define the Prandtl number 
\moneqstar  
 Pr =  {{ \rho \, \nu \, c_p}\over{\kappa}}  \, . 
\monendstar  
The conservation of mass, momentum and energy takes the form 
\moneq   \left\{ \begin{array} {l} 
\displaystyle  \partial_t \rho + \partial_x J =  0 \\ \displaystyle 
\partial_t J + \partial_x  ( \rho \, u^2 + p  ) 
- \partial_x  (  \rho \, \nu \, \partial_x u  )   = 0  \\ \displaystyle     
\partial_t  ( \rho \, E  ) + 
 \partial_x  (  \rho \, E \, u + p \, u  ) 
 - \partial_x  (  \rho \, \nu \, u \, \partial_x u  )  
- \partial_x  ( \kappa \, \partial_x T  )   =  0 \, . 
\end{array} \right.  \label{NS-1d}  \monend 
A first and natural idea proposed by  Alexander {\it et al.} \cite{ACS93},   
Khobalatte and Perthame \cite{KP94},  
McNamara and  Alder \cite{MNA93}, 
Shan and Chen \cite{SC93}  
is to use  two particle distributions. The first one ($f$) devoted to the 
conservation of  mass and momentum, the second one ($g$) 
 for the conservation of total energy. We consider  this idea by a coupling of  
fluid and  thermal D1Q3 lattice Boltzmann models. 
The conserved moments are those of the Navier-Stokes equations (\ref{NS-1d}):
\moneqstar   
\rho  =   \rho_0 \, \sum_j f_j  \,,\,\,  
 J =  \rho_0 \, \lambda \, \sum_j \, v_j \, f_j  \,,\,\,     
 \rho \, E =   \rho_0  \,  \lambda^2  \, \sum_j g_j  \, . 
\monendstar 
The nonconserved moments complete the set of moments for this double ``D1Q3Q3'':
\moneqstar   
   e  =  \rho_0  \,  \lambda^2 \,    ( f_+ + f_- - 2 \, f_0  )  \,,\,\,   
  \psi  = \rho_0  \,  \lambda^3 \, \sum_j \, v_j \, g_j \,,\,\,      
 \varepsilon   =  \rho_0  \,  \lambda^4  \,   ( g_+ + g_- - 2 \, g_0  ) \, . 
\monendstar 
The particle representation is now a vector with 6 components: 
\moneqstar   
  f_d  =   ( f_0 ,\, f_+ ,\, f_- ,\, g_0  ,\, g_+ ,\, g_-  ) ^{\rm  t}   
\monendstar 
and the vector of moments admits the expression  
\moneqstar   
    m  \equiv   ( \rho ,\, J ,\,  \rho \, E  ,\,  e  ,\, \psi   ,\, \varepsilon  )^{\rm  t}  \, . 
\monendstar 
We have $ m =  M \, f_d  $ with a matrix $M$ that is, up to a permutation,  the tensor 
product of the D1Q3 matrix with itself. 

\smallskip \noindent 
We focus now on previous works of  Guo {\it et al.} \cite{GZS07}, 
Li {\it et al.} \cite{LHWT07, LLHT12}, Karlin  {\it et al.} \cite{KSC13})
with the framework following strictly the approach of BKG \cite{BGK54}. 
These authors discretize with the  lattice Boltzmann method the set of 
kinetic equations 
\moneq 
\partial_t f_j  +  v_j \smb \nabla f_j  =  - {{1}\over{\tau_f}} \,   
\big( f_j - f_j^{\rm eq}  \big) \,,\quad 
\partial_t g_j +  v_j \smb \nabla g_j  =  - {{1}\over{\tau_g}} \,   
\big(  g_j - g_j^{\rm eq}  \big) 
+   {{Z}\over{\tau_f}} \,  \big( f_j - f_j^{\rm eq}  \big)   \, . 
 \label{bgk-couple}  \monend 
Remark that due to the $  {{Z}\over{\tau_f}} $ term in the right hand side of the 
second equation of (\ref{bgk-couple}), the relaxation step  is deeply transformed. 
We have adapted this idea for the multi relaxation times approach  \cite{Du16}. 
We distinguish three types of moments: 

\noindent 
(i) the conserved moments $W$ such that 
 $ W_k^*=  W_k  $ for  $ k < N$, 

\noindent 
(ii) the ``first non-conserved'' moments $m$ such that 
$ m_k^* =  (1-s_k) \, m_k \,+ \, s_k \,   m_k^{\rm eq}  $ 
and 

\noindent (iii) the ``second non-conserved''  moments such that 
\moneqstar   
 m_\ell^* =  m_\ell  \,-\, s_\ell \,  ( m_\ell - m_\ell^{\rm eq}  ) 
\,+\, K_{\ell k} \,  ( m_k -  m_k^{\rm eq}  )  \,, 
\monendstar 
 with  $  m_k $ in the first family of nonconserved moments. 
Our results \cite{Du16} are not entirely satisfying and we have changed our point of view.

\bigskip \bigskip   \noindent {\bf  6) \quad  
 Lattice Boltzmann algorithm for the simulation 

\quad $\,\,\,\,\,\,$ of the thermal Navier Stokes equations  }   

\noindent 
The important remark is that we can  re-interpret the second equation of  (\ref{bgk-couple}) as 
an usual advection of the particle distribution with a source term. 
%
In the following, we approach with a lattice Boltzmann scheme 
the full Navier Stokes equations as a system of conservative partial differential equations 
with a source term: 
\moneqstar 
 \partial_t W  +    \partial_x  F(W)  -   \partial_x  \big( \Phi(W, \nabla W) \big) =  S  \,. 
\monendstar 
This choice is motivated by the fact that the equivalent equations of a lattice Boltzmann scheme
are always conservative. 
In this contribution, we have chosen as ``conserved variables'' the  mass $\rho$, the momentum 
$ J \equiv  \rho  \, u $ and the volumic entropy.
The   volumic entropy $ \zeta $ is the product of the volumic mass $ \rho $ multiplied by the 
specific entropy  $  s $. Then  $ W = ( \rho ,\, J ,\, \zeta )  $. 
The  equation of state   
is  the one of a  polytropic perfect gas
\moneq 
 p = (\gamma - 1)  \, \rho \, i  \, =  \, \rho \, r \, T   
 \, = \, p_0 \, \Big( {{\rho}\over{\rho_0}} \Big)^\gamma \,\, 
\exp \Big( {{\gamma \, (s - s_0)}\over{c_p}} \Big)  \, . 
 \label{gaz-parfait}  \monend 
The sound velocity  $c$ satisfies 
\moneqstar 
 c^2   =  \gamma \, {{p}\over{\rho}}  =   \gamma \, (\gamma - 1)  \, i   
\monendstar 
as previously.   With this $(\rho,J,\zeta)$ formulation, 
the Navier Stokes equations express the conservation of mass and momentum, and the 
production of entropy:  
\moneq \left\{ \begin{array} {l} 
\displaystyle 
 \partial_t \rho + \partial_x J  = 0 \\ \displaystyle 
  \partial_t J + \partial_x \big( \rho \, u^2 + p \big) 
- \partial_x \big(  \rho \, \nu \, \partial_x u \big) = 0  \\ \displaystyle  
 \partial_t  \zeta  + 
 \partial_x \big( \zeta \, u  \big) 
 - \partial_x \Big(  {{\kappa}\over{T}} \, \partial_x T \Big)  = 
 {{\rho \, \nu}\over{T}} \, \big( \partial_x u \big)^2  
+   {{\kappa}\over{T^2}} \, \big( \partial_x T \big)^2 \, . 
 \end{array}  \right.  \label{NS-rho-J-zeta}  \monend   
We use two particle distributions $f$ and $g$ as previously. 
The first distribution $f$ devoted to  mass and momentum and the second distribution  $g$ 
to   volumic entropy. 
The conserved moments are defined from the double particle distribution according to 
\moneq  
\rho =  \rho_0 \, \big( f_0 + f_+ + f_- \big)   , \quad 
J = \rho_0 \, \lambda \,  \big( f_+ - f_- \big)  , \quad   
 \zeta   = \rho_0  \,  c_p \, \big( g_0 + g_+ + g_- \big)   \, .  
 \label{moments-conserves-rho-J-zeta}  \monend 
The nonconserved moments are defined thanks to the previous considerations: 
\moneq   
 e =   \rho_0  \,  \lambda^2 \,   \big( f_+ + f_- - 2 \, f_0 \big)   , \quad    
  \psi =  \rho_0  \,  c_p   \, \lambda \, ( g_+ - g_- )    , \quad    
 \varepsilon  =   \rho_0  \, c_p   \, \lambda^2  \,  \big( g_+ + g_- - 2 \, g_0 \big)  \, .  
 \label{moments-non-conserves-rho-J-zeta}  \monend 
We  apply the general theory of  multi relaxation times lattice Boltzmann schemes  
with the particle representation 
\moneqstar 
 f_d \equiv ( f_0 ,\,  f_+ ,\, f_- ,\, g_0  ,\, g_+ ,\, g_- ) ^{\rm  t}  
\monendstar 
and the momentum representation 
\moneqstar 
 m \equiv  ( \rho ,\, J ,\,  \rho \, E  ,\, e   ,\, \psi   ,\, \varepsilon )^{\rm  t}  \, . 
\monendstar 
We have  $ m  =  M_{\rm D1Q3Q3} \, f_d  $. 
The matrix $M_{\rm D1Q3Q3}$ between particles and moments
is given by the relation 
\moneq   
 M_{\rm D1Q3Q3} \,=\, \rho_0 \,\,  \left( \begin{array}{cccccc} 
\displaystyle   1 &   1 & 1 & 0 & 0 & 0 \cr
\displaystyle   0 &   \lambda & -\lambda & 0 & 0 & 0 \cr 
\displaystyle   0 & 0 & 0 & c_p  &  c_p  & c_p  \cr
\displaystyle   -2 \, \lambda^2  & \lambda^2 &  \lambda^2 & 0 & 0 & 0 \cr
\displaystyle   0 & 0 & 0 & 0 &   c_p \, \lambda  &  -c_p \, \lambda  \cr
\displaystyle   0 & 0 & 0 & -2 \, c_p \, \lambda^2  & c_p \, \lambda^2 &  c_p \lambda^2  
 \end{array} \right) \, . 
 \label{M-d1q3q3}  \monend 
We must now specify the 
equilibrium functions for non-conserved moments:
\moneq   
e^{\eq} =   e^{\eq}   ( \rho ,\, J ,\, \zeta   ) \,,\quad   
\psi^{\eq}  =   \psi^{\eq}   (  \rho ,\, J ,\,  \zeta ) \,,\quad     
\varepsilon^{\eq} =   \varepsilon^{\eq}    (  \rho ,\, J ,\,  \zeta  )  \, . 
 \label{equilibres-d1q3q3}  \monend 
The coefficients $ s_e $, $ s_\psi $ and $ s_\varepsilon $  for the relaxation of non-conserved moments, 
{\it i.e.} 
\moneq   
e^* = e +  s_e \,  ( e^{\eq} -  e  ) \,,\quad     
\psi^*  =  \psi  +  s_\psi  \,  ( \psi^{\eq} -  \psi  ) \,,\quad      
\varepsilon^*  =  \varepsilon  +  s_\varepsilon  \,  ( \varepsilon^{\eq} -  \varepsilon ) \, . 
 \label{relaxation-d1q3q3}  \monend 
The time iteration of the D1Q3Q3 scheme can be simply written: 
\moneq   
\left\{  \begin{array}{l} 
f_0 (x , \, t + \Delta t)  = f_0^* (x , \, t ) ,\,\,    
f_\pm (x , \, t + \Delta t)  = f_\pm^* (x \mp \Delta x , \, t )  \\  
g_0 (x , \, t + \Delta t)  = g_0^* (x , \, t )  ,\,\,      
g_\pm (x , \, t + \Delta t)  = g_\pm^* (x \mp \Delta x , \, t ) \, . 
 \end{array} \right. 
 \label{iteration-d1q3q3}  \monend 
The adaptation to the presence of source terms is described in 
\cite{DLT14}. The gradients  $ \partial_x u $ and $ \partial_x T $ 
are evaluated with second order centered finite differences.

\bigskip \bigskip   \noindent {\bf   7) \quad  
  Linearized Navier Stokes equations }   

\noindent 
An important step is the study of the linearized Navier Stokes equations. 
We consider a reference state 
 $ W_0 = ( \rho_0  , \, \rho \, u_0 ,\, \rho \, s_0 ) $, with 
the associated sound velocity $ c_0 $ satisfying 
$ c_0^2 = {{\gamma \, p_0}\over{\rho_0}} $.  
We linearize the Navier Stokes system (\ref{NS-rho-J-zeta}) and obtain without difficulty
\moneq  \left\{  \begin{array}{l} 
\displaystyle  \partial_t \rho + \partial_x J  = 0 \\ \displaystyle
 \partial_t J + \big( c_0^2 - u_0^2 - {{s_0 \, c_0^2}\over{c_p}}  \big) \, \partial_x  \rho 
+ 2 \, u_0 \, \partial_x J  +  {{c_0^2}\over{c_p}}   \, \partial_x  \zeta  
- \nu_0 \, \partial_x^2  \rho + \nu_0 \, u_0 \, \partial_x^2 J   =  0  \\ \displaystyle
 \partial_t  \zeta  - u_0 \, s_0 \, \partial_x  \rho + s_0  \, \partial_x J 
+ u_0   \, \partial_x  \zeta 
- {{\nu_0}\over{Pr}} \, \big( (\gamma-1) \, c_p - \gamma \, s_0 \big)  \, \partial_x^2  \rho
- \gamma \,  {{\nu_0}\over{Pr}} \,  \partial_x^2 \zeta  =  0  \, . 
 \end{array} \right.  \label{NS-rho-J-zeta-linearise}  \monend 
We first decouple and simplify the system (\ref{NS-rho-J-zeta-linearise}).
We obtain after this operation a simple system: 
\moneq  \left\{  \begin{array}{l} 
\displaystyle  \partial_t \rho + \partial_x J  = 0 \\  \displaystyle
 \partial_t J + 
( c_0^2 - u_0^2 ) \, \partial_x  \rho + 2 \, u_0 \, \partial_x J 
 - \nu_0 \, \partial_x^2  \rho + \nu_0 \, u_0 \, \partial_x^2 J   = 0    \\  \displaystyle
\partial_t  \zeta + u_0   \, \partial_x  \zeta 
- \gamma \,  {{\nu_0}\over{Pr}} \,  \partial_x^2 \zeta  \,= \, 0 \, . 
 \end{array} \right.  \label{NS-rho-J-zeta-linearise-decouple}  \monend 
The system of equations (\ref{NS-rho-J-zeta-linearise-decouple}) is nothing else than the 
juxtaposition of a linearized version of the system (\ref{edp-d1q3}) of fluid equations
and an advection-diffusion equation of the type (\ref{edp-thermic-ordre-2}). 
Then it is easy to derive an equation for the equilibrium moment: 
\moneqstar 
   {2\over3} \lambda^2 \, \rho  + {1\over3}   e^{\rm eq}  =   ( c_0^2 - u_0^2 ) \, \rho + 2 \, u_0 \, J \, .
\monendstar
The defect of conservation $ \theta_e $  can be evaluated in a pure algebraic way:
\moneqstar 
 \theta_e = 3 \, ( \lambda^2 - 3 \, u_0^2 -  c_0^2 )  \, \partial_x J 
- 6 \, u_0 \, (c_0^2  - u_0^2 ) \, \partial_x \rho \, \simeq \,  3 \, \lambda^2 \partial_x J   \, .
\monendstar
The  coefficient of relaxation $s_e$  is determined through the kinematic viscosity: 
$  \nu_0 = \sigma_e \, \lambda \, \Delta x  $. 
From the  scalar ``thermal'' equation, we deduce an expression for the equilibrium moment 
  $   \psi^{\rm eq}  = u_0 \, \zeta  $.
For the  defect of conservation  $ \theta_\psi $, we have 
\moneqstar 
\theta_\psi =  \Big(  {2\over3} \lambda^2  - u_0^2   \Big) \,  \partial_x \zeta 
+  {1\over3} \,  \partial_x \varepsilon^{\rm eq}     \, .
\monendstar
The compatibility at second order is provided  under the condition 
\moneq  
 \Big[ \Big(  {2\over3} \lambda^2  - u_0^2 \Big) \,  \zeta +  {1\over3} \,  \varepsilon^{\rm eq}  
\Big] \, \sigma_\psi \, \Delta t \,=\,  \gamma \,  {{\nu_0}\over{Pr}} \, \zeta  \, .  
 \label{compatibilite}  \monend 
Then the equilibrium of the moment $ \varepsilon $ is easy to determine:
\moneqstar 
\varepsilon^{\rm eq}  = 3 \Big(  {{\gamma}\over{Pr}}  \, {{\sigma_e}\over{\sigma_\psi}} -  {2\over3} 
+  {{u_0^2}\over{\lambda^2}} \Big) \, \lambda^2 \,  \zeta \, \equiv \, \alpha  \, \lambda^2 \,  \zeta  \, .
\monendstar     
In order to enforce stability $(-2 < \alpha < 1$), 
we suggest the following  link between two relaxations: 
\moneq   
\sigma_\psi =  {3\over2} \,  {{\gamma}\over{Pr}}  \, \sigma_e \, .  
 \label{lien-sigma-sigma}  \monend 

\smallskip \noindent 
We consider again the linearized Navier Stokes equations 
introduced in (\ref{NS-rho-J-zeta-linearise}). 
We maintain  the previous relations for relaxation coefficients: 
$ \nu_0 = \sigma_e  \, \lambda \, \Delta x $, $  \sigma_\varepsilon  = 1.5 $ 
and the relation (\ref{lien-sigma-sigma}). 
To assume compatibility between the equations (\ref{edp-ordre-2}) and 
 (\ref{NS-rho-J-zeta-linearise}), 
the equilibria $ e^{\rm eq} $, $  \psi^{\rm eq} $ and  $ \varepsilon^{\rm eq} $
are necessarily the following linear functions of the conserved variables  $ \rho $, $J$ and $ \zeta $: 
\moneq  \left\{ \begin{array}{l}   
\displaystyle e^{\rm eq}  = 
\Big( 3 \, \big( 1- {{s_0}\over{c_p}} \big) \, c_0^2 - 3 \, u_0^2 - 2 \, \lambda^2 \Big) \, \rho  
+ 6  \, u_0  \, J  +  3  \, {{u_0^2}\over{c_p}} \, \zeta  \\ \displaystyle  
 \psi^{\rm eq}   =  - u_0 \, s_0 \, \rho  + s_0 \, J  +   u_0   \, \zeta  \\ \displaystyle
  \varepsilon^{\rm eq}   =  
\Big( -3 \, (s_0\, c_0)^2 - 6 \, s_0 \, u_0^2  + 3 \, s_0 \, c_0^2 + 2 - 2 \, s_0 - {{2}\over{\gamma}} \Big) 
\, \, \rho  \\ \displaystyle \qquad \qquad  \qquad  \qquad \qquad  \qquad
+  \, 6 \, u_0 \, s_0  \, J  +   3 \,  \big( c_p \,  u_0^2 +  s_0 \, c_0^2  \big) \, \zeta  \, . 
 \end{array} \right. 
 \label{equilibres-NS-linearise}  \monend 
For the simulation of a  simple linear wave with the following parameters: 
\moneqstar \left\{ \begin{array}{l} 
 \gamma = 1.4 \,, \,\, Pr = 1 \,, \,\,  c_0 =  {{\lambda}\over{2}}  \,, \,\,   
 u_0 = 0  \,, \,\,    s_0 = 0  \,, \\ 
s_e = 1.9  \,, \,\,  
\nu =  6.579 \, 10^{-4}  \,, \,\,   \Delta x =   {{1}\over{40}}  \,, \,\,  T_f =  120 \, \Delta t \,, 
 \end{array} \right. \monendstar 
we have observed numerical stability.
This fundamental property was also realized with the important modification of the parameters:
$ u_0 =  0.15 \, \lambda  $ and $ s_0 = 0.2 \, c_p  $.

\bigskip \bigskip   \noindent {\bf     8) \quad  
 D1Q3Q3 lattice Boltzmann scheme for the volumic entropy 

\quad $\,\,\,\,\,\,$  Navier Stokes equations  }   

\noindent 
The matrix $M_{\rm D1Q3Q3}$ between particles and moments
is still given by the relation (\ref{M-d1q3q3}). 
The equilibrium of non-conserved moments (\ref{equilibres-d1q3q3}) 
is parameterized by nonlinear functions. 
We must also  specify 
the coefficients $ s_e $, $ s_\psi $ and $ s_\varepsilon $  for the relaxation (\ref{relaxation-d1q3q3}) 
of the non-conserved moments. The discrete time iteration of the D1Q3Q3 scheme 
follows the relations (\ref{iteration-d1q3q3}). 
The Taylor expansion method at first order proposes  partial differential equations  
satisfied by the conserved variables 
$  W = ( \rho ,\, J ,\, \zeta \equiv \rho \, s ) $: 
\moneqstar 
\partial_t W_k + \Lambda_{k  \ell} \, \partial_x  m_\ell^{\eq} = {\rm O}(\Delta t) \monendstar 
with   a momentum-velocity tensor $\Lambda $ given in (\ref{Lambda-general}). 
We have for this D1Q3Q3 lattice Boltzmann scheme 
\moneq  
\Lambda  =    \left( \begin{array}{cccccc}
   0 &   1 & 0 & 0 & 0 & 0 \cr
   {{2 \, \lambda^2}\over{3}}  &  0 & 0 &   {{1}\over{3}}  & 0 & 0 \cr 
   0 & 0 & 0 & 0  &  1 & 0 \cr
   0  & \lambda^2 &  0  & 0 & 0 & 0 \cr
   0 & 0 &  {{2 \, \lambda^2}\over{3}}  & 0 & 0  &   {{1}\over{3}} \cr
   0 & 0 & 0 & 0  & \lambda^2 &  0
\end{array} \right)  \, . 
 \label{Lambda-d1q3q3}  \monend 
The equivalent equations at first order  take the form 
\moneq  
\partial_t \rho + \partial_x J =  {\rm O}(\Delta x) \,,\quad  
\partial_t J + \partial_x \Big(  {{2 \, \lambda^2}\over{3}} \, \rho 
+  {{1}\over{3}} \, e^{\eq} \Big)  =  {\rm O}(\Delta x) \,,\quad   
 \partial_t \zeta  + \partial_x  \psi^{\eq}  =  {\rm O}(\Delta x) \,. 
 \label{edp1-d1q3q3}  \monend 
They are compared  
 to the Navier stokes equations of gas dynamics  (\ref{NS-rho-J-zeta}). 
By identification of first order terms, the ``energy'' at equilibrium $  e^{\eq} $
is given by the relation (\ref{energie-equilibre})
and  
\moneq  
  \psi^{\eq} =  \zeta  \, u\,. 
 \label{psi-equilibre}  \monend 
%
Thus two equilibria for nonequilibrium moments are fixed.

\smallskip \noindent 
For the second order analysis, the H\'enon's coefficients \cite{He87} 
are defined according to 
\moneq  
\sigma_{e} =  {{1}\over{s_{e}}} -   {1\over2} \,, \quad 
\sigma_\psi =  {{1}\over{s_\psi}} -   {1\over2} \,, \quad  
\sigma_\varepsilon =  {{1}\over{s_\varepsilon}} -   {1\over2} \, . 
 \label{henon-d1q3q3}  \monend 
The defects of conservation (\ref{thetas})
have to be estimated before the determination 
of the second order equivalent partial  differential equations (\ref{edp-ordre-2}). 
%
For the present  D1Q3Q3 lattice Boltzmann scheme, we have  
\moneq  \left\{ \begin{array}{l}   
\displaystyle \partial_t \rho + \partial_x J =  {\rm O}(\Delta t) \\ \displaystyle
\partial_t J + \partial_x \big( \rho \, u^2 + p \big)  \,= \, {{\Delta t}\over{3}}  \, 
\sigma_{e} \,\, \partial_x \theta_{e} + {\rm O}(\Delta x^2)  \\ \displaystyle
  \partial_t \zeta +  \partial_x \big( \zeta \, u  \big) \,= \, 
\Delta t  \, \, \sigma_\psi \,\, \partial_x \theta_\psi  + {\rm O}(\Delta x^2) \, . 
\end{array} \right. \label{edp2-d1q3q3}  \monend 
The first defect of conservation 
\moneqstar 
 \theta_e  \equiv   \partial_t e^{\eq} + \lambda^2 \, \partial_x J 
\monendstar 
can be approximated by 
$  3 \, \lambda^2 \, \rho \, \partial_x u  $.  Then the relation 
\moneq   
\nu = \sigma_e \, \lambda \, \Delta x
 \label{nu-d1q3q3}  \monend 
determines the relaxation coefficient $ s_e $ thanks to  (\ref{henon-d1q3q3}).
For the second  defect of conservation 
\moneqstar 
\theta_\psi  \equiv 
\partial_t \psi^{\eq} +  {2\over3} \, \lambda^2 \, \partial_x \zeta  + 
  {{1}\over{3}} \, \partial_x \varepsilon^{\eq} 
\monendstar 
 we have,  after integration  by parts 
and using the first order equations  (\ref{edp1-d1q3q3}),  
\moneqstar 
 \theta_\psi  = p \, \partial_x s + \partial_x \big[  {1\over3} \, 
\big( 2 \, \lambda^2 \, \zeta +  \varepsilon^{\eq} \big)
-  ( ( \rho \, u^2 + p ) \, s ) \big] \, . 
\monendstar 
This expression is approximated by 
$  \partial_x  \big[  {1\over3} \, \big( 2 \, \lambda^2 \, \zeta +  
\varepsilon^{\eq} \big) - \big( ( \rho \, u^2 + p ) \, s \big) \big]   $.  
The entropy equation (third equation of (\ref{NS-rho-J-zeta})) 
is compared with the third equation of (\ref{edp2-d1q3q3}): 
\moneq   
 \partial_t \zeta +  \partial_x \big( \zeta \, u \big) 
-  \partial_x^2  \big\{ \sigma_\psi \, \Delta t \, 
\big[  {1\over3} \, \big( 2 \, \lambda^2 \, \zeta +  
\varepsilon^{\eq} \big) - \big( ( \rho \, u^2 + p ) \, s \big) \big]  \big\}   
=  {\rm O}(\Delta x^2) \, . 
 \label{eq-entropie-d1q3q3}  \monend 
Then 
\moneqstar     
\sigma_\psi \, \Delta t \, 
\Big[  {1\over3} \, \big( 2 \, \lambda^2 \, \zeta +  
\varepsilon^{\eq} \big) - \big( ( \rho \, u^2 + p ) \, s \big) \Big] = 
       \kappa  \, {\rm log } \Big( {{T}\over{T_0}} \Big)  \, . 
\monendstar    
With the relations  (\ref{lien-sigma-sigma}) and (\ref{nu-d1q3q3}), we obtain 
for the third  momentum for entropy  
\moneq     
\varepsilon^{\eq} = {{2 \, \rho \, c_p \, \lambda^2}\over{\gamma}} \, 
{\rm log } \Big( {{T}\over{T_0}} \Big) + 
\big[ 3 \, \big( \rho \, u^2 + p \big) - 2 \, \lambda^2 \, \rho \big]  \, s \, . 
 \label{eq-energie-entropie-d1q3q3}  \monend 
The relation (\ref{eq-energie-entropie-d1q3q3}) is clearly  non-trivial. 
It is an interesting property of the present scheme to  have 
determined  this algebraic relation  in the full  nonlinear framework.
For the polytropic perfect gas, the pressure is given by (\ref{gaz-parfait}) 
and the relation (\ref{eq-energie-entropie-d1q3q3}) becomes 
\moneqstar     
\varepsilon^{\eq}  = 2 \,  \lambda^2 \big[ \rho \, (s - s_0) 
+ \big( 1 - {{1}\over{\gamma}} \big) \, c_p \, \rho \, 
{\rm log } \big( {{\rho}\over{\rho_0}} \big) \big] +  
\big[ 3 \, \big( \rho \, u^2 + p \big) - 2 \, \lambda^2 \, \rho \big]   \, s \, .  
\monendstar    
%
The present D1Q3Q3 algorithm can be implemented in the following way:

\noindent
{\bf (i)} consider the particle distribution 
$ f_d \equiv ( f_0 ,\,  f_+ ,\, f_- ,\, g_0  ,\, g_+ ,\, g_- ) ^{\rm  t} $, 

\noindent
{\bf (ii)} compute the six moments 
$ m \equiv  ( \rho ,\, J ,\,  \zeta  ,\, e   ,\, \psi   ,\, \varepsilon )^{\rm  t}  $,  
by  $ m = M_{D1Q3Q3} \,  f_d $ with $ M_{D1Q3Q3} $ defined by 
the relation (\ref{M-d1q3q3}), 

\noindent
{\bf (iii)} fix the relaxation coefficients $s_e$, $s_\psi$, $s_\varepsilon$   of the non-conserved moments  
 $ e $, $ \psi $, $ \varepsilon $ through the associated $\sigma's$ thanks to 
(\ref{henon-d1q3q3}) with $ \sigma_e$ determined by (\ref{nu-d1q3q3}), $ \sigma_\psi $
according to (\ref{lien-sigma-sigma}) and $ s_\varepsilon = 1.5$, 

\noindent
{\bf (iv)} determine the moments   $ e^{\rm eq} $, $ \psi^{\rm eq} $, $ \varepsilon^{\rm eq} $ at equilibrium
with the relations (\ref{energie-equilibre}), (\ref{psi-equilibre}) and 
(\ref{eq-energie-entropie-d1q3q3}) respectively, 

\noindent
{\bf (v)} compute the moments   $ e^* $, $ \psi^* $, $ \varepsilon^* $ after relaxation 
 with the relations (\ref{relaxation-d1q3q3}), 

\noindent
{\bf (vi)} compute the particle distribution $ f_d^*$ after relaxation thanks to 
\moneqstar     
 f_d^* =  (M_{D1Q3Q3})^{-1}  \,   ( \rho ,\, J ,\,  \zeta  ,\, e^*   ,\, \psi^*   ,\, \varepsilon^* )^{\rm  t} 
\,,  \monendstar 
{\bf (vii)} evaluate the gradients  $ \partial_x u $ and $ \partial_x T $ 
with second order centered finite differences, add a source term to the momentum  $\zeta$
as proposed by the third equation  in (\ref{NS-rho-J-zeta}) and the algorithm described in \cite{DLT14},

\noindent 
{\bf (viii)} iterate the scheme in time according to  (\ref{iteration-d1q3q3}). 

%
\bigskip 
\noindent {\bf     9) \quad  
 First numerical experiments }   

\smallskip \smallskip \smallskip \smallskip 
 \centerline  {\includegraphics[width=.70 \textwidth]{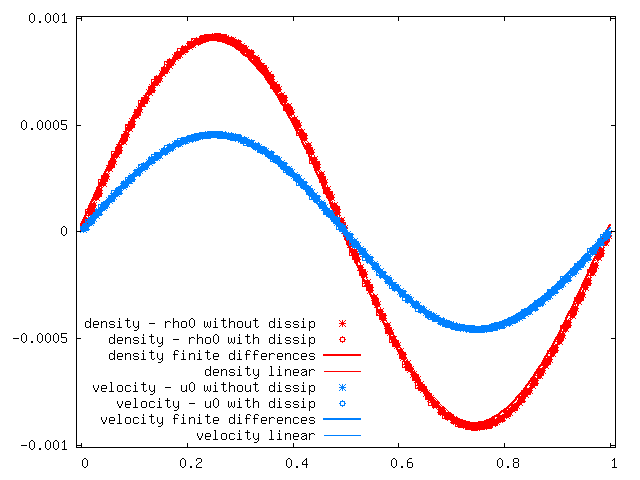}}

 \centerline {\includegraphics[width=.70 \textwidth]{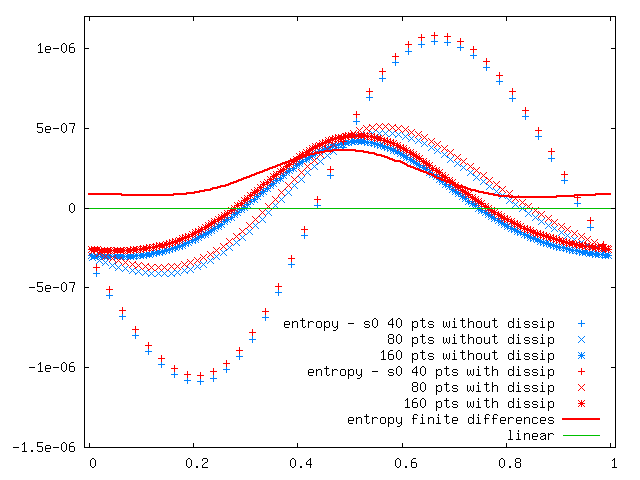}} 


\noindent  
{\small {\bf Figure 1.} 
Progressive linear wave, periodic boundary conditions, $ \delta \rho = 0.001 \, \rho_0 $.
Density and velocity (top). Entropy (bottom). }

\smallskip \smallskip \smallskip \smallskip 

\centerline  {\includegraphics[width=.70 \textwidth]{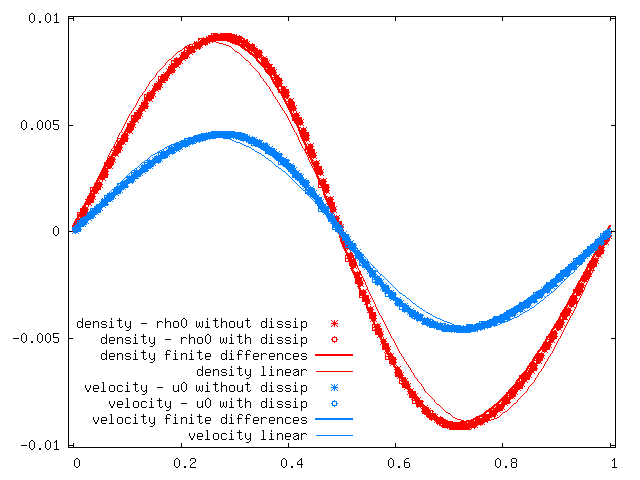}}

\centerline {\includegraphics[width=.70 \textwidth]{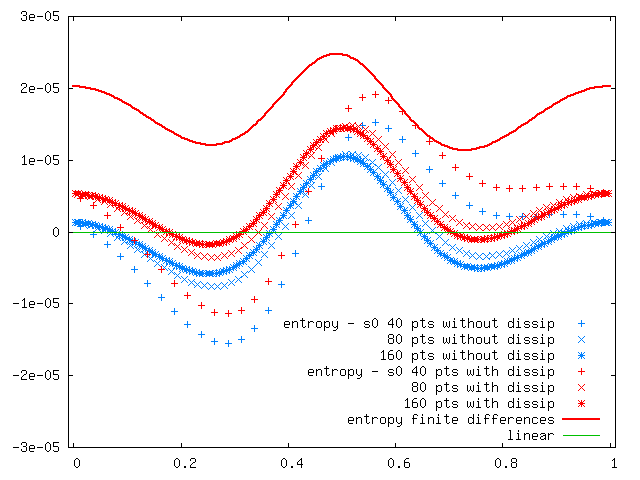}}  

\noindent 
{\small {\bf Figure 2.}  
Progressive nonlinear wave, periodic boundary conditions, $ \delta \rho = 0.01 \, \rho_0 $.
Density and velocity (top). Entropy (bottom). }  

\bigskip \noindent 
We have used the algorithm described in the previous section for the simulation of periodic waves. 
The general parameters are 
\moneqstar 
\gamma = 1.4 \,,\,\, Pr = 1 \,,\,\,  c_0 =  {{\lambda}\over{2}} \,,\,\,    
  s_0 = 0 \,,\,\,   T_{\rm max} =  3 \,,\,\,  \nu =  6.579 \, 10^{-4} \, .  
\monendstar 
We have used three  meshes with 
$   \Delta x =  {{1}\over{40}} $, $  {{1}\over{80}} $,  $  {{1}\over{160}} $. 
We present three numerical experiments. The initial condition is a simple acoustic wave 
with an initial variation $ \delta \rho $ of density.
For the first experiment (see Figure~1), 
the result is essentially described 
by the linearized equations. The variation of entropy is infinitesimal. 
%

\smallskip \smallskip \smallskip \smallskip 

  \centerline {\includegraphics[width=.70 \textwidth]{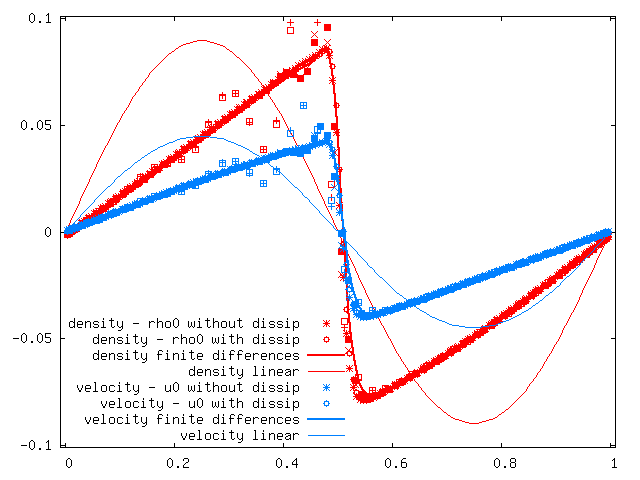}}

  \centerline {\includegraphics[width=.70 \textwidth]{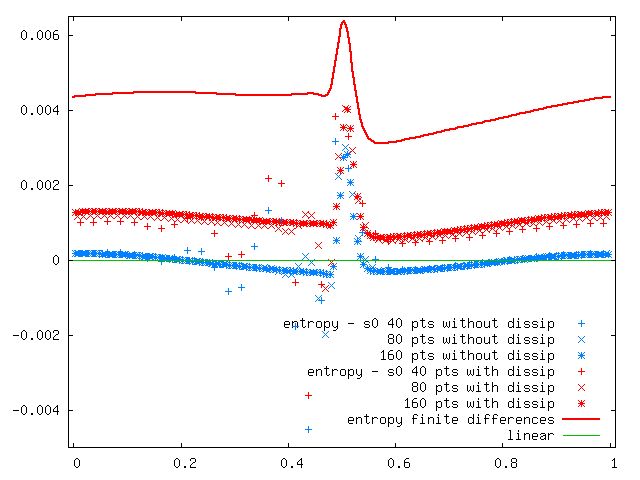}} 

\noindent 
{\small {\bf Figure 3.}  
Strong nonlinear wave, periodic boundary conditions, $ \delta \rho = 0.1 \, \rho_0 $.
Density and velocity (top). Entropy (bottom). }   

\smallskip  \smallskip \noindent 
With the second experiment described in Figure~2, 
nonlinear effects are clearly visible
for the density and velocity  results.   
In order to validate these effects, 
we have developed a   finite difference software 
based on the formulation (\ref{NS-rho-J-zeta}) 
and using  second order centered finite differences and  explicit first order time integration.   
The  variation of entropy is very small and is of good quality.  
With $ \,  \delta \rho = {{\rho_0}\over{100}} $, a shock wave is generated 
(see Figure~3).  
The production of  entropy is essentially 
localized in the region of high variation of the  fields. 
On Figure~4,  
we present the transient 
 evolution of specific entropy and 
 specific energy. Without the dissipation source term of (\ref{NS-rho-J-zeta}), 
the total energy if not perfectly conserved and 
the specific entropy is quasi constant. 
With this source term, the defect for total energy is very small and the production of entropy 
remains moderate in comparison  to our  finite difference simulator.
The extension of the previous work to two and  three space dimensions 
is the next step of this work, typically with the coupled lattice Boltzmann schemes 
D2Q9-D2Q5 and  D3Q19-D3Q7. 

 \smallskip \smallskip \smallskip \smallskip 
 \centerline{  {\includegraphics[width=.70 \textwidth]{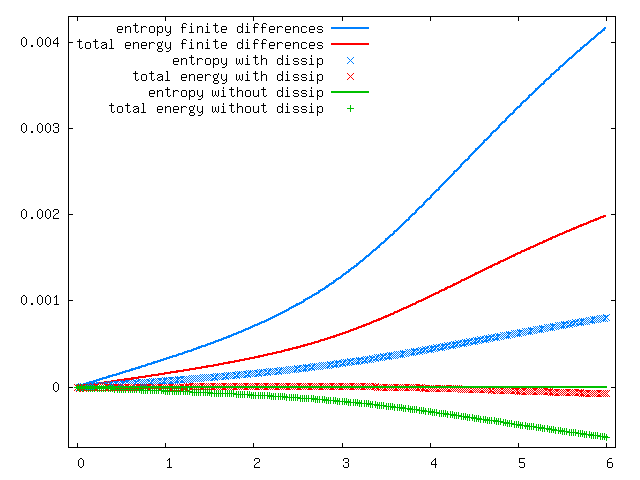}} }

\noindent 
{\small {\bf Figure 4.} 
Strong nonlinear wave ($ \delta \rho = 0.1 \, \rho_0 $). Time evolution of entropy 
and total energy.} 

\bigskip \bigskip      \noindent {\bf    Acknowledgements  }   

\noindent 
This work is supported  by the French ``Climb'' Oseo project. 
A part of this work has been realized during the stay of two of us at the 
Beijing  Computational Science Research Center. We thank the colleagues of CSRC 
for their hospitality. Last but not least, the authors  thank the referee  who suggested several points 
 in need of improvement.


\bigskip \bigskip      \noindent {\bf    References }   


\end{document}